\documentclass[12pt]{article}
\usepackage{amsmath,amssymb}
\usepackage{graphicx}
 
\begin{document}

\title{Global Regularity of the 3D Axi-symmetric Navier-Stokes Equations
with Anisotropic Data}
\author{Thomas Y. Hou\thanks{Applied and Comput. Math, Caltech, Pasadena,
CA 91125. Email: hou@acm.caltech.edu.}
\and Zhen Lei \thanks{Applied and Comput. Math, Caltech, Pasadena,
CA 91125. Email: zhenlei@acm.caltech.edu.}
\and Congming Li \thanks{Department of
Applied Mathematics, University of Colorado, Boulder, CO. 80309. Email: cli@colorado.edu}}

\maketitle

\begin{abstract}
In this paper, we study the 3D axisymmetric Navier-Stokes Equations with
swirl. We prove the global regularity of the 3D Navier-Stokes
equations for a family of large anisotropic initial data. Moreover,
we obtain a global bound of the solution in terms of its initial data
in some $L^p$ norm. Our results also reveal some interesting
dynamic growth behavior of the solution due to the interaction
between the angular velocity and the angular vorticity fields.
\end{abstract}

\section{Introduction.}

Despite a great deal of effort by many mathematicians and physicists,
the question of whether the solution of the 3D Navier-Stokes equations
can develop a finite time singularity from a smooth initial condition
with finite energy remains one of the most outstanding open problems
\cite{Clay}.  A main difficulty in obtaining the global regularity of
the 3D Navier-Stokes equations is due to the presence of the vortex
stretching, which is absent for the 2D problem. So far,
the global regularity of the 3D Navier-Stokes equations has
been obtained only for initial data which are small in some
scaling invariant norm \cite{Lady70,Temam01,BM02}. But the analysis
for small initial data does not generalize to the 3D Navier-Stokes
with large data.
A more refined analysis which takes into account the special
nature of the nonlinearities and the anisotropic nature of
the solution near the region of blow-up seems to be needed.

In this paper, we study the global regularity of the axisymmetric
Navier-Stokes equations with large initial data that have
anisotropic scaling.
Let $u^\theta$ and $\omega^\theta$ be the angular velocity and
vorticity components of the 3D axisymmetric Navier-Stokes equations.
We consider initial data for $u^\theta$ and $\omega^\theta$
that have the following scaling property:
\begin{eqnarray}
u^\theta(r,z,0) = \frac{1}{\epsilon^{1-\delta}} U_0
(\epsilon r,z ), \quad
\omega^\theta(r,z,0)& = & \frac{1}{\epsilon^{1-\delta}} W_0
(\epsilon r,z ),
\label{Data-1}
\end{eqnarray}
%\begin{eqnarray}
%u^\theta(r,z,0) = \frac{1}{\epsilon} U_0\left (\frac{r}{\epsilon^{1-\delta}},
%\frac{z}{\epsilon} \right ), \quad
%\omega^\theta(r,z,0) = \frac{1}{\epsilon^2} W_0
%\left (\frac{r}{\epsilon^{1-\delta}},\frac{z}{\epsilon}\right ),
%\label{Data-2}
%\end{eqnarray}
where $r=\sqrt{x^2+y^2}$, $\delta$ and $\epsilon $
are some small positive parameters,
and the rescaled profiles $U_0$ and $W_0$ are bounded in $L^{2p}$
and $L^{2q}$ respectively for some $p$ and $q$ with $p=2q$.
We note that these initial data are not small. In fact, we have
\begin{eqnarray*}
\| {\bf u}_0 \|_{L^2(\mathbb{R}^2\times [0,1])} \|\nabla {\bf
u}_0\|_{L^2(\mathbb{R}^2\times [0,1])} =
\frac{C_0}{\epsilon^{4-2\delta}} \gg 1,
\end{eqnarray*}
for $\epsilon$ small, where ${\bf u}_0$ is the initial velocity vector
(we use bold letters to denote vector fields throughout this paper).
Thus the classical regularity analysis for small initial data does
not apply to these sets of anisotropic initial data.

In this paper, we prove the global regularity of the 3D axisymmetric
Navier-Stokes equations for both initial data given by (\ref{Data-1})
and (\ref{Data-2}) by exploring the anisotropic structure of the solution
for $\epsilon$ small.  We also obtain a global bound
on $\|u^\theta\|_{L^{2p}}$ and $\|\omega^\theta\|_{L^{2q}}$ in
terms of their initial data. Note that by using the scaling
invariance property of the Navier-Stokes equations, our global
regularity result also applies to the following rescaled
initial data
\begin{eqnarray}
u^\theta(r,z,0) = \frac{1}{\epsilon^{2-\delta}} U_0
\left (r,\frac{z}{\epsilon} \right ), \quad
\omega^\theta(r,z,0)& = & \frac{1}{\epsilon^{3-\delta}} W_0
\left (r, \frac{z}{\epsilon} \right ),
\label{Data-2}
\end{eqnarray}
and
\begin{eqnarray}
u^\theta(r,z,0) = \frac{1}{\epsilon} U_0\left (\frac{r}{\epsilon^{1-\delta}},
\frac{z}{\epsilon} \right ), \quad
\omega^\theta(r,z,0) = \frac{1}{\epsilon^2} W_0
\left (\frac{r}{\epsilon^{1-\delta}},\frac{z}{\epsilon}\right ).
\label{Data-3}
\end{eqnarray}
For the rescaled initial data (\ref{Data-2}) and (\ref{Data-3}),
our analysis suggests
that $\|u^\theta\|_{L^{2p}}$ may experience a rapid growth dynamically
due to the contribution from $\omega_0^\theta$. On the other hand,
the growth of $\|\omega^\theta\|_{L^{2q}}$ is more moderate.
See Remark 2 and (\ref{bound-u1})-(\ref{bound-w1}) in Section 3 for
more discussions.
This phenomenon is similar to the solution behavior of the 1D model
derived by Hou-Li in \cite{HL06b} where they show that $\omega^\theta$
is essentially bounded by its initial data while $u^\theta$ can
experience large growth dynamically due to the contribution from
$(u_0^\theta)_z$ and $\omega_0^\theta$.

Note that the parameters $\epsilon $ in the
initial data (\ref{Data-1})-(\ref{Data-2}) and $\delta$ in (\ref{Data-3})
measure the degree of anisotropy of the initial data.
If $\delta = 0$, then the initial data (\ref{Data-3})
become isotropic, i.e.
\[
{\bf u}_0(x,y,z) = \frac{1}{\epsilon} {\bf U}_0\left (\frac{x}{\epsilon},
\frac{y}{\epsilon},\frac{z}{\epsilon} \right ).
\]
Our analysis breaks down when there is no anisotropic scaling in
the initial data, i.e. $\delta = 0$. Clearly, if the analysis
could be extended to the case of $\delta = 0$, one would prove the
global regularity of the 3D axisymmetric Navier-Stokes equations
for general initial data by using the scaling invariance
property of the Navier-Stokes equations. It is interesting
to note that by using an anisotropic scaling of the initial
data, we turn the global regularity of the 3D Navier-Stokes
equations into a critical case of $\delta=0$.

We would like to emphasize that our global regularity results are
obtained on a regular size domain, $\mathbb{R}^2\times [0,1]$, for
initial data (\ref{Data-1}). In this sense, our results are
different from those global regularity results obtained for a thin
domain, $\Omega_\epsilon = Q_1 \times [0,\epsilon]$ with $Q_1$
being a bounded domain in $\mathbb{R}^2$. We remark that the
global regularity of the 3D Navier-Stokes equations in a thin
domain of the form $\Omega_\epsilon$ has been studied by Raugel
and Sell in a series of papers \cite{RS93a,RS94,RS93b}. They prove
the global regularity of the 3D Navier-Stokes equations under the
assumption that $\| \nabla {\bf u}_0\|_{L^2(\Omega_\epsilon)}^2
\le C_0 \ln \frac{1}{\epsilon}.$ This is an improvement over the
classical global regularity result for small data, which requires
$\| \nabla {\bf u}_0\|_{L^2(\Omega_\epsilon)}^2 \le C^* \epsilon$
\cite{RS93a}. We may interpret our global regularity result with
initial data (\ref{Data-2}) as a result on a generalized thin
domain. Note that the initial data given by (\ref{Data-2}) satisfy
the following bound: $\| \nabla {\bf
u}_0\|_{L^2(\Omega_\epsilon)}^2 = C_0\epsilon^{-5+2\delta}$ (here
$\delta>0$ can be made arbitrarily small), which is much larger
than the corresponding bound $ C_0 \ln \frac{1}{\epsilon}$
required by the global regularity analysis of Raugel and Sell in
\cite{RS93a,RS94,RS93b}.

This paper is motivated by the desire to understand how the local
anisotropic structure of the solution near the region of a potential
singularity may lead to the depletion of vortex stretching, thus
preventing the formation of a finite time singularity. To fully
characterize the local solution structure of a potential singularity
would be an extremely difficult task. Such study would shed useful
light into our understanding of the dynamic depletion of vortex
stretching and the global regularity of the Navier-Stokes equations.
There has been some encouraging progress along this
direction recently. In particular, Necas, Ruzicka and Sverak
\cite{NRS96} and Tsai \cite{Tsai98} have ruled out the possibility
of an isotropic self-similar blow-up of the 3D Navier-Stokes equations.
More recently,  Chen-Strain-Tsai-Yau \cite{Tsai07} prove that
for axisymmetric 3D Navier-Stokes equations, if the velocity field
${\bf u}$ satisfies a scaling invariant blow-up rate
\[
|{\bf u} (x,t)| \leq C_*/\sqrt{r^2 + (T-t)},
\]
where $r=\sqrt{x^2+y^2}$, then such blow-up is not possible
at $t=T$. This to some extent excludes the locally self-similar
isotropic blow-up of the axisymmetric Navier-Stokes equations.
A related result with a stronger assumption can be found in
\cite{HL07} for the 3D Navier-Stokes equations.

There have been also some numerical evidences on the anisotropic scaling
of the solution near the region of a potential blow-up of the
3D Euler and Navier-Stokes equations, see e.g.
\cite{Kerr93,BP94,Pelz97,GP00,Pelz01,Kerr05,HL06,HL07b}. In particular,
the recent work of Kerr \cite{Kerr05} and Hou-Li \cite{HL06,HL07b}
gives some detailed description of the anisotropic scaling of the vorticity
field near the inner core region of the maximum vorticity for
two slightly perturbed antiparallel vortex tubes. The
motivation for considering the rescaled initial data
(\ref{Data-2}) and (\ref{Data-3}) is to understand whether
such locally anisotropic blow-up is possible for the Navier-Stokes
equations.

We remark that Hou-Li have recently proved the global regularity
of the 3D axisymmetric Navier-Stokes equations with some
large anisotropic initial data \cite{HL06b}. The results obtained
in this paper complement those presented in \cite{HL06b}. For
more discussions on the axisymmetric solutions of the Navier-Stokes
equations, we refer to \cite{CL02} and \cite{Tsai07} and the references
cited there. The 2D Boussinesq equations are closely related to the
3D axisymmetric Navier-Stokes equations with swirl (away from the
symmetry axis). Recently, Hou-Li \cite{HL05} and Chae \cite{Chae05}
have proved independently the global existence of the 2D viscous
Boussinesq equations with viscosity entering only in the fluid
equation, but the density equation remains inviscid.

The remaining part of the paper is organized as follows.
We give the formulation of the problem in Section 2. In
Section 3, we prove a useful property of Riesz operators in
$\mathbb{R}^4 \times \bf{T}^1$, where $\bf{T}^1$ is the
one-dimensional torus with periodicity 1. This property is
needed to prove our global regularity analysis of the axisymmetric
Navier-Stokes equations with anisotropic initial data.
Finally, we will present and prove the main result of this paper
in Section 4.

\section{Formulation}

Consider
the 3D axi-symmetric incompressible Navier-Stokes equations with swirl.
\begin{equation}
  \left\{ \begin{array}{l}
   {\bf u}_t + ({\bf u} \cdot \nabla ) {\bf u} = - \nabla p + \Delta {\bf u},\\
    \nabla \cdot {\bf u} = 0,\\
    {\bf u}|_{t = 0} = {\bf u}_0 ( {\bf x}), \quad {\bf x} = (x,y,z).
  \end{array} \label{nse} \right.
\end{equation}

Let
\[
{\bf e}_r = \left ( \frac{x}{r}, \frac{y}{r}, 0 \right ),
\;\;
{\bf e}_\theta = \left (-\frac{y}{r}, \frac{x}{r}, 0 \right ),\;\;
{\bf e}_z = \left ( 0, 0, 1 \right ),
\]
be three unit vectors along the radial,
the angular, and the $z$ directions respectively,
$r=\sqrt{x^2+y^2}$. We will decompose the
velocity field as follows:
\begin{equation}
{\bf u} = v^r(r,z,t) {\bf e}_r + u^\theta (r,z,t) {\bf e}_\theta +
v^z (r,z,t) {\bf e}_z.
\end{equation}
In the above expression, $u^\theta$ is called the {\it swirl} component
of the velocity field ${\vec u}$.
The vorticity field can be expressed similarly
\begin{equation}
{\boldsymbol{\omega}} = -(u^\theta )_z (r,z,t) {\bf e}_r + \omega^\theta (r,z,t)
{\bf e}_\theta +
\frac{1}{r} (r u^\theta )_r (r,z,t) {\bf e}_z,
\end{equation}
where $\omega^\theta = v_z^r - v^z_r $.

One can derive evolution equations for $u^\theta$ and $\omega^\theta$ as follows
(see e.g. \cite{BM02,CL02}).
\begin{eqnarray}
&& u^\theta_t + v^r u^\theta_r + v^z u^\theta_z =  \left ( \nabla^2 - \frac{1}{r^2} \right )
u^\theta  - \frac{1}{r} v^r u^\theta ,\label{eqn-u} \\
&& \omega^\theta_t + v^r \omega^\theta_r + v^z \omega^\theta_z
 =  \left ( \nabla^2 - \frac{1}{r^2} \right )
\omega^\theta +\frac{1}{r} \left ( (u^\theta)^2 \right )_z + \frac{1}{r}
v^r \omega^\theta ,\label{eqn-w} \\
&& - \left ( \nabla^2 - \frac{1}{r^2} \right ) \psi^\theta = \omega^\theta ,
\label{eqn-psi}
\end{eqnarray}
where $\psi^\theta$ is the angular component of the stream function,
$v^r$ and $v^z$ can be expressed in terms of the angular stream function
$\psi^\theta$ as follows:
\begin{equation}
v^r = - \frac{\partial \psi^\theta}{\partial z}, \quad
v^z =  \frac{1}{r} \frac{\partial }{\partial r} (r \psi^\theta ),
\end{equation}
and
$\nabla^2$ is defined as
\begin{equation}
\nabla^2 = \partial_r^2 + \frac{1}{r} \partial_r + \partial_z^2 .
\label{laplacian}
\end{equation}
Note that equations (\ref{eqn-u})-(\ref{eqn-psi}) completely
determine the evolution of the 3D axisymmetric Navier-Stokes
equations once the initial condition is given.

By the well-known Caffarelli-Kohn-Nirenberg theory \cite{CKN82},
the singularity set of any suitable weak solution of the 3D Navier-Stokes
equations has one-dimensional Hausdorff measure zero.
Thus, in the case of axisymmetric 3D Navier-Stokes
equations with swirl, if there is any singularity, it must be along the
symmetry axis, i.e. the $z$-axis. Therefore, we should focus our effort
to understand the possible singular behavior of the 3D Navier-Stokes
equations near the symmetry axis at $r=0$.

As observed by Liu and Wang in \cite{LW04}, any smooth solution
of the 3D axisymmetric Navier-Stokes equations must satisfy
the following compatibility condition at $r=0$:
\begin{equation}
u^\theta(0,z,t)=\omega^\theta(0,z,t)=\psi^\theta(0,z,t)=0.
\label{eqn-bc}
\end{equation}
Thus we can rewrite $u$, $\omega$ and $\psi$ as follows:
\begin{equation}
\label{u1-eqn}
u^\theta (r,z,t)= r u_1(r,z,t), \quad \omega^\theta (r,z,t) = r \omega_1(r,z,t),
\quad \psi^\theta (r,z,t) = r \psi_1(r,z,t).
\end{equation}
In \cite{HL06b}, Hou-Li have derived
the following equivalent system for $(u_1, \omega_1,\psi_1)$:
\begin{eqnarray}
&& ( u_1)_t + v^r (u_1)_r + v^z (u_1)_z = 2 \psi_{1z} u_1 +
(u_{1zz} + u_{1rr} + \frac{3 u_{1r}}{r} )
\label{u1-dyn} \\
&& ( \omega_1)_t + v^r (\omega_1)_r + v^z (\omega_1)_z = (( u_1)^2)_z +
(\omega_{1zz} + \omega_{1rr} + \frac{3 \omega_{1r}}{r} )
\label{w1-dyn} \\
&& - ( \psi_{1zz} + \psi_{1rr} + \frac{3 \psi_{1r}}{r} )
= \omega_1 ,
\label{psi1-dyn}
\end{eqnarray}
where $v^r$ and $v^z$ are defined as
\begin{equation}
v^r = - \frac{\partial ( r \psi_1)}{\partial z}, \quad
v^z =  \frac{1}{r} \frac{\partial }{\partial r} (r^2 \psi_1 ).
\label{eqn-v3r}
\end{equation}
We note that the incompressibility condition implies that
\begin{equation}
(r v^r)_r + (r v^z)_z = 0,
\label{incompress}
\end{equation}
which can be verified directly from (\ref{eqn-v3r}).

\section{A Useful Property of Riesz operators on $\mathbb{R}^4 \times \bf{T}^1$}

Before we present our global regularity analysis, we would like to
state and prove two technical lemmas regarding the property of a
Riesz operator on $\mathbb{R}^4 \times \bf{T}^1$, where
$\mathbb{T}^1$ is the one-dimensional torus with periodicity 1.
We first recall the following weighted Calderon-Zygmund inequality
for a singular integral operator with a weight function
which is in the $A_p$ class (see Stein \cite{Stein93} pp. 194-217
for details). Let $K$ be a Riesz operator in
$\mathbb{R}^n$ and $w(x)$ be a weight in the $A_p$ class
(see page 194 of \cite{Stein93} for definition). One can extend
the Calderon-Zygmund inequality for the singular integral
operator with the integral having weight function $w(x)$.
Specifically, for $1 < p < \infty$, there exists a uniform
constant $C$ such that
\begin{equation}\label{s2}
\int_{\mathbb{R}^n}|K*f|^pw(x)dx \leq
C^p\int_{\mathbb{R}^n}|f|^pw(x)dx,
\end{equation}
for all $f \in L^p(\mathbb{R}^n)$.

In this section, we prove a similar property for a Riesz operator
on $\mathbb{R}^4 \times \bf{T}^1$. Let $\Delta = \Sigma_{i =
1}^4\partial_{x_i}^2 + \partial_z^2$. Assume that $f(r,z)$ belongs to
the weighted $L^p(\mathbb{R}^4 \times \mathbb{T}^1)$ space with a
weight function $w(x) = \frac{1}{r^2}$. Let $u = u(r, z)$, $r =
\sqrt{x_1^2 + x_2^2 + x_3^2 + x_4^2}$, be the solution of
\begin{equation}\label{s3}
- \Delta u = f,
\end{equation}
with periodic boundary conditions along the $z$ direction with
period 1. We will prove the following lemma:

\vspace{0.1in} \noindent
 {\bf Lemma 1.} {\it Assume that $w(x) = w(r)$ belongs to
the $A_p$ class in $\mathbb{R}^5$. Then there exists
 a uniform positive constant $C > 0$ such that
\begin{equation}\label{s4}
\int_{\mathbb{R}^{n - 1} \times
\bf{T}^1}|\nabla^2u|^pw(x)dx_1\cdots dx_{n - 1}dz \leq
C\int_{\mathbb{R}^{n - 1} \times \bf{T}^1}|f|^pw(x)dx_1\cdots
dx_{n - 1}dz,
\end{equation}
provided that $f$ and $u$ are both $L^p$ integrable with the weight
$w(x)$ for $1 < p < \infty$.}

\vspace{0.1in} \noindent {\bf Proof of Lemma 1.} Let $\lambda > 0$
and denote the one-dimensional torus with periodicity
$\frac{1}{\lambda}$ by $\frac{1}{\lambda}\mathbb{T}^1$. For $(x_1,
\cdots, x_4, z) \in \mathbb{R}^4 \times
\frac{1}{\lambda}\mathbb{T}^1$, define
\begin{equation}\nonumber
\begin{cases}
u^\lambda(x_1, \cdots, x_4, z) =  u(\lambda x_1, \cdots,
\lambda x_4, \lambda z),\\
f^\lambda(x_1, \cdots, x_4, z) =  \lambda^2f(\lambda x_1, \cdots,
\lambda x_4, \lambda z).
\end{cases}
\end{equation}
It is easy to see that
\begin{equation}\label{s10}
- \Delta u^\lambda = f^\lambda,\quad (x_1, \cdots, x_4, z) \in
\mathbb{R}^4 \times \frac{1}{\lambda}\mathbb{T}^1.
\end{equation}
We claim that there exists a uniform constants $C_1$ and $C_2$
independent of $\lambda$ such that
\begin{equation}\label{s11}
\|(\nabla^2u^\lambda) w(x)^{\frac{1}{p}}\|_{L^p(\mathbb{R}^4
\times \frac{1}{\lambda}\mathbb{T}^1)} \leq C_1\|f^\lambda
w(x)^{\frac{1}{p}}\|_{L^p(\mathbb{R}^4 \times
\frac{1}{\lambda}\mathbb{T}^1)} + C_2\|u^\lambda
w(x)^{\frac{1}{p}}\|_{L^p(\mathbb{R}^4 \times
\frac{1}{\lambda}\mathbb{T}^1)}.
\end{equation}
To prove this, we introduce a smooth cutoff function $\phi(z)$
which satisfies $0 \leq \phi(z) \leq 1$, $\phi(z) = 0$ for
$z \leq - 1$ or $z \geq 2$, and $\phi(z) = 1$ for $0 \leq z \leq 1$.
For $\lambda \in (0, 1]$,
we define $\phi^\lambda(z) = \phi(\lambda z)$. A simple computation
gives
\begin{equation}\nonumber
- \Delta (\phi^\lambda u^\lambda) = \phi^\lambda f^\lambda -
2\nabla \phi^\lambda \cdot \nabla u^\lambda -
u^\lambda\Delta\phi^\lambda,\quad (x_1, \cdots, x_4, z) \in
\mathbb{R}^{5}.
\end{equation}
Since $w(x)$ belongs to the $A_p$ class, the weighted
Calderon-Zygmund inequality in $\mathbb{R}^n$ (see (\ref{s2})) implies
that
\begin{eqnarray}\nonumber
&&\|(\nabla^2u^\lambda) w(x)^{\frac{1}{p}}\|_{L^p(\mathbb{R}^4
  \times\frac{1}{\lambda}\mathbb{T}^1)} \leq
  \|\nabla^2(\phi^\lambda u^\lambda) w(x)^{\frac{1}{p}}\|_{L^p(\mathbb{R}^5)}\\\nonumber
&&\leq C\big\|\big(\phi^\lambda f^\lambda -
  2\nabla\phi^\lambda \cdot \nabla u^\lambda
 - u^\lambda\Delta\phi^\lambda\big)w(x)^{\frac{1}{p}}\big\|_{L^p(\mathbb{R}^5)}\\\nonumber
&&\leq C\big(\|f^\lambda w(x)^{\frac{1}{p}}\|_{L^p(\mathbb{R}^4
  \times \frac{1}{\lambda}\mathbb{T}^1)} + \lambda^2\|
  u^\lambda w(x)^{\frac{1}{p}}\|_{L^p(\mathbb{R}^4 \times
  \frac{1}{\lambda}\mathbb{T}^1)}\big)\\
\label{Weighted}
&&\quad +\ \frac{1}{2}\|(\nabla^2u^\lambda)
w(x)^{\frac{1}{p}}\|_{L^p(\mathbb{R}^4 \times
  \frac{1}{\lambda}\mathbb{T}^1)},
\end{eqnarray}
where we have used the following estimate
(note that $\phi =\phi(z)$, $w(x)=w(r)$ and $u(x)=u(r,z)$):
\begin{eqnarray}\nonumber
&&\big\|\nabla\phi^\lambda \cdot \nabla u^\lambda
 w(x)^{\frac{1}{p}}\big\|_{L^p(\mathbb{R}^5)}
 \leq C\lambda\Big(\int_{\mathbb{R}^4\times \frac{1}{\lambda}\mathbb{T}^1}
|\partial_z u^\lambda|^pw(x)
 dx\Big)^{\frac{1}{p}}\\\nonumber
&&\leq C\lambda\Big(
\int_{\mathbb{R}} w(r) r^3dr
\int_{\frac{1}{\lambda}\mathbb{T}^1}
|\partial_z u^\lambda|^p dz
  \Big)^{\frac{1}{p}} \\\nonumber
&&\leq C_p\lambda^2\|u^\lambda w(x)^{\frac{1}{p}}
  \|_{L^p(\mathbb{R}^4 \times \frac{1}{\lambda}\mathbb{T}^1)}
  + \frac{1}{4}\|(\partial_z^2u^\lambda)
  w(x)^{\frac{1}{p}}\|_{L^p(\mathbb{R}^4 \times
  \frac{1}{\lambda}\mathbb{T}^1)}.
\end{eqnarray}
In the last step of the above estimate, we have used the
one-dimensional Sobolev inequality along the $z$-direction:
\[
\| u_z\|_{L^p( \frac{1}{\lambda}\mathbb{T}^1)}^p
\leq \frac{C_p \lambda}{C} \|u\|_{L^p( \frac{1}{\lambda}\mathbb{T}^1)}^p
+ \frac{1}{4C\lambda} \|\partial_z^2 u\|_{L^p( \frac{1}{\lambda}\mathbb{T}^1)}^p,
\]
for some positive constant $C_p$ depending on $p$.
Thus, one has
\begin{eqnarray}\label{s13}
\|(\nabla^2u^\lambda) w(x)^{\frac{1}{p}}\|_{L^p(\mathbb{R}^4
  \times\frac{1}{\lambda}\mathbb{T}^1)} \leq C\big(\|f^\lambda w(x)^{\frac{1}{p}}\|_{L^p(\mathbb{R}^5)} + \lambda^2\|
  u^\lambda w(x)^{\frac{1}{p}}\|_{L^p(\mathbb{R}^4 \times
  \frac{1}{\lambda}\mathbb{T}^1)}\big)
\end{eqnarray}
for some uniform positive constant $C > 0$. For $\lambda
> 1$, we denote $N$ by the biggest integer less that $\lambda$ and
compute
\begin{equation}\nonumber
- \Delta (\phi u^\lambda) = \phi f^\lambda - 2\nabla\phi \cdot
\nabla u^\lambda - u^\lambda\Delta\phi,\quad (x_1, \cdots, x_4, z)
\in \mathbb{R}^{5}.
\end{equation}
Similarly, by using the weighted Calderon-Zygmund inequality in
$\mathbb{R}^5$, we obtain
\begin{eqnarray}\nonumber
&&N\|(\nabla^2u^\lambda)
w(x)^{\frac{1}{p}}\|_{L^p(\mathbb{R}^4\times
  \frac{1}{\lambda}\mathbb{T}^1)} \leq C
  \|\nabla^2(\phi u^\lambda) w(x)^{\frac{1}{p}}\|_{L^p(\mathbb{R}^5)}\\\nonumber
&&\leq C\big\|\big(\phi f^\lambda - 2\nabla\phi \cdot \nabla
  u^\lambda - u^\lambda\Delta\phi\big)w(x)^{\frac{1}{p}}\big\|_{L^p(\mathbb{R}^5)}\\\nonumber
&&\leq CN\big(\|f^\lambda w(x)^{\frac{1}{p}}\|_{L^p(\mathbb{R}^4
  \times\frac{1}{\lambda}\mathbb{T}^1)} + \|
  u^\lambda w(x)^{\frac{1}{p}}\|_{L^p(\mathbb{R}^4 \times
  \frac{1}{\lambda}\mathbb{T}^1)}\big)\\\nonumber
&&\quad  +\ \frac{N}{2}\|(\nabla^2u^\lambda)
w(x)^{\frac{1}{p}}\|_{L^p(\mathbb{R}^4 \times
  \frac{1}{\lambda}\mathbb{T}^1)},
\end{eqnarray}
which gives
\begin{eqnarray}\label{s15}
\|(\nabla^2u^\lambda) w(x)^{\frac{1}{p}}\|_{L^p(\mathbb{R}^4
  \times
  \frac{1}{\lambda}\mathbb{T}^1)} \leq C\big(\|f^\lambda w(x)^{\frac{1}{p}}
  \|_{L^p(\mathbb{R}^4 \times
  \frac{1}{\lambda}\mathbb{T}^1)} + \|
  u^\lambda w(x)^{\frac{1}{p}}\|_{L^p(\mathbb{R}^4\times
  \frac{1}{\lambda}\mathbb{T}^1)}\big).
\end{eqnarray}
By \eqref{s13} and \eqref{s15}, we prove the claim \eqref{s11}.

Note that
\begin{equation}\nonumber
\begin{cases}
\|(\nabla^2u^\lambda) w(x)^{\frac{1}{p}}\|_{L^p(\mathbb{R}^4\times
  \frac{1}{\lambda}\mathbb{T}^1)} =
  \lambda^{2 - \frac{5}{p}}\|(\nabla^2u) w(x)^{\frac{1}{p}}\|_{L^p(\mathbb{R}^4\times
  \mathbb{T}^1)},\\
\|f^\lambda w(x)^{\frac{1}{p}}\|_{L^p(\mathbb{R}^4\times
  \frac{1}{\lambda}\mathbb{T}^1)} =
  \lambda^{2 - \frac{5}{p}}\|f w(x)^{\frac{1}{p}}\|_{L^p(\mathbb{R}^4\times \mathbb{T}^1)},\\
\|u^\lambda w(x)^{\frac{1}{p}}\|_{L^p(\mathbb{R}^4\times
  \frac{1}{\lambda}\mathbb{T}^1)} =
  \lambda^{- \frac{5}{p}}\|u w(x)^{\frac{1}{p}}\|_{L^p(\mathbb{R}^4\times
  \mathbb{T}^1)}.
\end{cases}
\end{equation}
Thus, we deduce from \eqref{s11} that
\begin{equation}\label{s16}
\|(\nabla^2u) w(x)^{\frac{1}{p}}\|_{L^p(\mathbb{R}^4 \times
\mathbb{T}^1)} \leq C_1\|f w(x)^{\frac{1}{p}}\|_{L^p(\mathbb{R}^4
\times \mathbb{T}^1)} + \frac{C_2}{\lambda^2}\|u
w(x)^{\frac{1}{p}}\|_{L^p(\mathbb{R}^4 \times \mathbb{T}^1)}.
\end{equation}
By letting $\lambda \rightarrow \infty$ in \eqref{s16}, we prove
the lemma.

\vspace{0.1in} \noindent

Now we prove the following Lemma, which plays an important
role in our global regularity analysis:

\vspace{0.1in}
\noindent {\bf Lemma 2.} {\it Assume that
$\omega_1$ and $\psi_1$ are in
$L^p(\mathbb{R}^{2} \times \mathbb{T}^1)$ with $ 1 < p < \infty$,
and $\psi_1$ is the
solution of \eqref{psi1-dyn} with periodic boundary condition along
the $z$-direction with period 1. Then we have
\begin{equation}
\| \psi_{1zz} \|_{L^p(\mathbb{R}^{2} \times \mathbb{T}^1)}
\leq C \| \omega_1 \|_{L^p(\mathbb{R}^{2} \times \mathbb{T}^1)},
\label{Imbeding}
\end{equation}
for $ 1 < p < \infty$.}

\vspace{0.1in}
\noindent {\bf Proof of Lemma 2.} To prove
(\ref{Imbeding}), we observe that
\[
\Delta = \frac{\partial^2}{\partial z^2} +
\frac{\partial^2}{\partial r^2} + 3 \frac{\partial}{\partial r}
\]
is an axisymmetric Laplacian operator in five space dimensions,
where $r = \sqrt{x_1^2+x_2^2+x_3^3+x_4^2}$ for
$x=(x_1,x_2,x_3,x_4,z)\in R^5$. Let $n=5$.
Since $-\Delta \psi_1 = \omega_1$, Lemma 1 implies that
\begin{equation}
\int_{\mathbb{R}^4 \times
\bf{T}^1}|\partial_z^2\psi_1|^pw(x)dx_1\cdots dx_4dz \leq
C\int_{\mathbb{R}^4\times \bf{T}^1}|\omega_1|^pw(x)dx_1\cdots
dx_4dz. \label{CZL}
\end{equation}
provided that $w(x)= \frac{1}{r^2}$ belongs to the $A_p$ class.
We will show that $w$ belongs to
the $A_p$ class. In terms of the cylindrical coordinate, the
above inequality (\ref{CZL}) can be written as
\begin{equation}
\left ( \int_0^1\int_0^\infty |\psi_{1zz}|^p r^3r^{-2} dr dz \right )^{\frac{1}{p}}
\leq C \left ( \int_0^1\int_0^\infty |\omega_1 |^p r^3r^{-2} dr dz \right
)^{\frac{1}{p}}, \label{CZL1}
\end{equation}
which is exactly \eqref{Imbeding}.

It remains to show that $w(x) = \frac{1}{r^2}$ belongs to the
$A_p$ class. In the following, we will prove that
$w(x) = \frac{1}{r^\alpha}$
belongs to the $A_p$ class for $ 0 < \alpha < 4$ (see page 194 of
\cite{Stein93} for definition). To see this, let $B_R$ be a ball
centered at $x=0$ with radius $R$. Denote by $|B_R|$ the volume of
$B_R$, and $p'$ be the conjugate of $p$, i.e. $\frac{1}{p}
+\frac{1}{p'} = 1$. We have
\begin{eqnarray*}
&& \left ( \frac{1}{|B_R|} \int_{B_R} w(x) dx \right )
\left ( \frac{1}{|B_R|} \int_{B_R} w(x)^{-\frac{p'}{p}} dx \right )^{\frac{p}{p'}}\\
&& \leq \frac{C}{R^5} \left ( \int_0^R \int_0^{\sqrt{R^2-z^2}}
r^{3-\alpha} drdz \right ) \left ( \frac{1}{R^5} \int_0^R dz
\int_0^{\sqrt{R^2-z^2}} r^{3+\alpha\frac{p'}{p}} dr dz
\right )^{\frac{p}{p'}} \\
&& \leq \frac{C}{R^5} \left ( \int_0^R
(R^2-z^2)^{\frac{4-\alpha}{2}} dz \right ) \left ( \frac{1}{R^5}
\int_0^R {(R^2-z^2)}^{2+\frac{\alpha p'}{2p}} dz
\right )^{\frac{p}{p'}} \\
&& \leq \frac{C}{R^5} R^{5-\alpha} \left ( \frac{1}{R^5}
R^{5+\alpha\frac{p'}{p}} \right )^{\frac{p}{p'}} \leq C R^{\alpha
- \alpha} = C,
\end{eqnarray*}
independent of $R$ for $0 < \alpha < 4$. Similarly we can prove
the same estimate for $B_R$ centered at any point. This completes
the proof of Lemma 2.

\section{The main result and its proof}

Now we state our main result in this paper.

\vspace{0.1in}
\noindent
{\bf Theorem 1.} {\it Let $0 < \delta < 1$ be a given number which
could be arbitrarily small. Assume that
\begin{eqnarray}
u^\theta(r,z,0) = \frac{1}{\epsilon^{1-\delta}} U_0
(\epsilon r,z ), \quad
\omega^\theta(r,z,0) =  \frac{1}{\epsilon^{1-\delta}} W_0
(\epsilon r, z),
\label{UW-init1}
\end{eqnarray}
where $U_0(r,\tilde{z}) = r U_1(r,\tilde{z})$,
$W_0(r,\tilde{z}) = r W_1(r,\tilde{z})$, and
$U_0,\; U_1 \in L^{2p} (R^2\times [0,1])$,
$W_0, \; W_1 \in L^{2q} (R^2\times[0,1])$ with $p=2q$ and $q> 1/\delta$.
Further, we assume that the initial velocity field
${\bf u}_0 \in L^2(R^2\times [0,1])$ and the initial conditions for
$\psi^\theta$, $u^\theta$ and $\omega^\theta$ are odd and periodic in
$z$ with period 1. Then there exists $\epsilon_0 (\delta, U_0, W_0) >0$
such that for all $0<\epsilon \leq \epsilon_0$ the 3D axisymmetric Navier-Stokes
equations with initial data given by (\ref{UW-init1})
have a unique global regular solution. }

\vspace{0.1in}
\noindent
{\bf Remark 1.}
By using Theorem 1, we can easily obtain the global regularity
of the initial data (\ref{Data-2}) by using a scaling argument.
Specifically, if we denote by $u$ the solution of the Navier-Stokes
equations with initial data (\ref{UW-init1}), then
$u_\epsilon ({\bf x},t) \equiv \frac{1}{\epsilon}
u\left (\frac{\bf x}{\epsilon},\frac{t}{\epsilon^2}\right )$
is the solution of the Navier-Stokes equations with the
rescaled initial data (\ref{Data-2}). Thus the global
regularity of the Navier-Stokes
equations with initial data (\ref{Data-2}) follows immediately
from Theorem 1. Using a slightly different
rescaling, we can prove the global regularity of the
Navier-Stokes equations with initial data (\ref{Data-3}).
See Remark 4 at the end of this section for more detail.

\vspace{0.1in}
\noindent
{\bf Remark 2.}
It is interesting to derive the corresponding global bound
for the initial data (\ref{Data-2}).
Let $(u^\theta, \;\omega^\theta)$ be the solution corresponding
to the initial data (\ref{UW-init1}), and
$(u_\epsilon^\theta, \;\omega_\epsilon^\theta)$ be the
solution corresponding to the initial data (\ref{Data-2}).
If we denote by $u_1^\epsilon = u_\epsilon^\theta/r$,
$\omega_1^\epsilon = \omega_\epsilon^\theta/r$, and
$u_1=u^\theta/r$, $\omega_1 = \omega^\theta/r$, then we have
\[
u_1^\epsilon (r,z,t) = \frac{1}{\epsilon^2} u_1 (\frac{r}{\epsilon},\frac{z}{\epsilon},
\frac{t}{\epsilon^2}),
\quad
\omega_1^\epsilon (r,z,t) = \frac{1}{\epsilon^3} \omega_1
(\frac{r}{\epsilon},\frac{z}{\epsilon}, \frac{t}{\epsilon^2}).
\]
Using the above scaling relationship and
substituting $f= u_1^p$ and $g = \omega_1^q$ into
the global estimate (\ref{global-bound}), we obtain the following
global estimate for $u_1^\epsilon$ and
$\omega_1^\epsilon$ with the initial data (\ref{Data-2}):
\begin{eqnarray}
\| u_1^\epsilon (t) \|_{L^{2p}(\Omega_\epsilon)}
& \leq & \epsilon^{-2+\delta + 1/2p}\| U_1 \|_{L^{2p}(\Omega_\epsilon)}
+ (2/C_q)^{\frac{1}{2p}} \epsilon^{-2+\delta/2 + 1/2p}
\| W_1 \|_{L^{2q}(\Omega_\epsilon)}^{\frac{1}{2}} \label{bound-u1}, \\
\| \omega_1^\epsilon (t) \|_{L^{2q}(\Omega_\epsilon)}
& \leq & \epsilon^{-3+\delta + 1/2q}\| W_1 \|_{L^{2q}(\Omega_\epsilon)}
+ (C_q/2)^{\frac{1}{2q}} \epsilon^{-3+2\delta + 1/2q}
\| U_1 \|_{L^{2p}(\Omega_\epsilon)}^2\label{bound-w1},
\end{eqnarray}
where $\Omega_\epsilon = R^2\times [0,\epsilon]$.
Note that the first term on the right hand side of (\ref{bound-u1})
is equal to $\| u_1^\epsilon (0) \|_{L^{2p}(\Omega_\epsilon)}$ and
is smaller than the second term by a factor of $O(\epsilon^{\delta/2})$.
This implies that $\|u_\epsilon^\theta\|_{L^{2p}}$ may experience
a rapid growth dynamically for $\epsilon$ small due to the
contribution from the initial condition of
$\omega_\epsilon^\theta$. Similarly, the first term on the right
hand side of (\ref{bound-w1}) is equal to
$\| \omega_1^\epsilon (0) \|_{L^{2p}(\Omega_\epsilon)}$, but is
larger than the second term by a factor of $O(\epsilon^{-\delta})$.
Thus the growth of $\|\omega_\epsilon^\theta\|_{L^{2q}}$ is more moderate.
This phenomenon is similar to the
solution behavior of the 1D model derived by Hou-Li in \cite{HL06b}
where they show that $(u_{1z})^2+\omega_1^2$ satisfies a
maximum principle. However, $u_1$ can experience large
growth dynamically due to the contribution from
$u_{1z}$ and $\omega_1$ at $t=0$.

\vspace{0.1in} \noindent {\bf Remark 3.} It is worthwhile to point out
that all the functions involved are in fact smooth for $t>0$ based
on our assumptions and the dynamic control that we obtain. Here, the
key is the dynamic control of $ \| u_1 (t) \|_{L^{2p}} \leq C$ and $
\| w_1 (t) \|_{L^{p}} \leq C$. With these two estimates, the
solution can be smoothly extended to all times. More specifically,
the work of \cite{CKN82} and \cite{Lin98} implies that the
velocity field $\textbf{u}(x,t)$ is square integrable and smooth in
the region where $r>1$. Thus $\textbf{u}(x,t)$ has a bounded $L^p$
norm in the region $r>1$ for $2 \leq p \leq \infty$. On the other hand,
we have $|u^\theta | \leq C|u_1|$ and $|w^{\theta}| \leq C |w_1|$
in the region where $r$ is bounded. Therefore,
our dynamic control of $ \| u_1 (t) \|_{L^{2p}} \leq C$ (note that
$2p \geq 3$ by the assumption of Theorem 1) and
$ \| w_1 (t) \|_{L^{p}} \leq C$ lead to the following
{\it a priori} bound on the gradient of $u^r{\bf e}_r + u^z{\bf e}_z$
(see \cite{NRS96}):
\begin{eqnarray}\nonumber
&&\|\nabla (u^r{\bf e}_r + u^z{\bf e}_z)\|_{L^{p}(\mathbb{R}^2 \times
  \mathbb{T}^1, r < 1.5)}\\\nonumber
&&\leq C\big(\|\nabla\times (u^r{\bf e}_r + u^z{\bf e}_z)
  \|_{L^{p}(\mathbb{R}^2 \times \mathbb{T}^1, r < 2)} +
  \|\nabla\cdot (u^r{\bf e}_r + u^z{\bf e}_z)\|_{L^{p}(\mathbb{R}^2 \times
  \mathbb{T}^1, r < 2)}\\\nonumber
&&\quad +\ \|u^r{\bf e}_r + u^z{\bf e}_z\|_{L^2(\mathbb{R}^2 \times
  \mathbb{T}^1,  r < 2)}\big) \\\nonumber
&&\leq C\big(\|\omega^\theta {\bf e}_\theta
  \|_{L^{p}(\mathbb{R}^2 \times \mathbb{T}^1, r < 2)}
  + \|u^r{\bf e}_r + u^z{\bf e}_z\|_{L^2(\mathbb{R}^2 \times
  \mathbb{T}^1,  r < 2)}\big)\\\nonumber
&&< \infty,
\end{eqnarray}
where we have used the fact that
$ \omega^\theta {\bf e}_\theta = \nabla \times (u^r{\bf e}_r + u^z{\bf e}_z)$
and $\quad \nabla \cdot (u^r{\bf e}_r + u^z{\bf e}_z) = 0$.
The Sobolev embedding theory and the fact that $u^r{\bf e}_r + u^z{\bf e}_z \in
L^2(\mathbb{R}^2 \times  \mathbb{T}^1,  r < 2)$ imply that
$u^r{\bf e}_r + u^z{\bf e}_z \in L^3(\mathbb{R}^2 \times \mathbb{T}^1,r < 1.5)$ since
$p\geq \frac{3}{2}$.
Then the well-known regularity result (see
\cite{ISS03}) implies that $\textbf{u}(t)$ is smooth
and $L^p$ integrable with $p \geq 2$ for $t>0$.

\vspace{0.1in}
\noindent
{\bf Proof of Theorem 1.}
By the assumptions of Theorem 1, we know that the initial velocity field
${\bf u}_0 \in L^2(R^2\times [0,1])$ and that
$u^\theta|_{t=0} \in L^{2p}(R^2\times [0,1])$,
$\omega^\theta|_{t=0} \in L^{2q}(R^2\times [0,1])$ with
$p=2q$ and $q > 1/\delta$. Using these properties of the
initial data, one can show that there exists a finite time $T>0$
such that the axisymmetric Navier-Stokes has a unique regular
solution for $0< t \leq T$. In the rest of the proof, we will
perform a dynamic estimate for the regular solution for
$0< t \leq T$ and show that $T$ can be made arbitrarily large.

Denote $\Omega = R^2\times [0,1]$. Define
$f=|u_1|^p$ and $g=|\omega_1|^q$, with $p=2q$. Multiplying
(\ref{u1-dyn}) by $|u_1|^{2p-2}u_1$ and integrating over $\Omega$, we
get after using (\ref{incompress}) that
\begin{eqnarray}
\frac{1}{2p} \frac{d}{dt} \int_{\Omega} f^2 r dr dz =
\int_{\Omega} \psi_{1z} |u_1|^{2p} r dr dz + \int_{\Omega}
|u_1|^{2p-2}u_1 ( u_{1zz} + u_{1rr} + \frac{3 u_{1r}}{r} ) r dr dz .
\end{eqnarray}
Note that
\begin{eqnarray*}
\int_{\Omega} |u_1|^{2p-2}u_1 u_{1zz} r drdz & =& - (2p-1) \int_{\Omega}
|u_1|^{2p-2} (|u_1|_z)^2 r dr dz \\
& = & - (2p-1) \int_{\Omega} \left ( |u_1|^{p-1} |u_1|_z \right )^2 r dr dz \\
& = & - \frac{(2p-1)}{p^2} \int_{\Omega} (f_z)^2 rdr dz .
\end{eqnarray*}
Similarly, we obtain
\begin{eqnarray*}
\int_{\Omega} |u_1|^{2p-2}u_1 \frac{(r u_{1r})_r}{r} r drdz & =&
\int_{\Omega} |u_1|^{2p-2}u_1 (r u_{1r})_r dr dz \\
& = & - (2p-1) \int_{\Omega} |u_1|^{2p-2} (|u_1|_r)^2 r dr dz \\
& = & - \frac{(2p-1)}{p^2} \int_{\Omega} (f_r)^2 rdr dz ,
\end{eqnarray*}
and
\begin{eqnarray*}
\int_{\Omega} |u_1|^{2p-2}u_1 \frac{2 u_{1r}}{r} r drdz & =&
2 \int_{\Omega} |u_1|^{2p-2} u_1 u_{1r} dr dz  =
\frac{1}{p} \int_{\Omega} \left (|u_1|^{2p} \right )_r dr dz \\ \\
& = & - \frac{1}{p} \int_0^1 |u_1|^{2p} (0,z,t) dz
= - \frac{1}{p} \int_0^1 f^2(0,z,t)^2 dz .
\end{eqnarray*}
Therefore, we have
\begin{eqnarray}
\int_{\Omega} |u_1|^{2p-2} u_1( u_{1zz} + u_{1rr} + \frac{3 u_{1r}}{r} ) r dr dz
= - \frac{(2p-1)}{p^2} \int_{\Omega} |\nabla f |^2 r dr dz -
\frac{1}{p} \int_0^1 f^2(0,z,t)^2 dz .
\end{eqnarray}

On the other hand, if we let $\tilde p = 2q$, $\tilde q = \tilde p/(\tilde p -1)$,
we get
\begin{eqnarray*}
\int_{\Omega} \psi_{1z} f^2 r dr dz &\leq & \| \psi_{1z}\|_{L^{\tilde p}(\Omega)}
\| f^2 \|_{L^{\tilde q}(\Omega)} \leq
\| \psi_{1zz} \|_{L^{\tilde p}(\Omega)}
\| f^2 \|_{L^{\tilde q} (\Omega)} \\
&\leq &  C
\| \omega_1 \|_{L^{\tilde p}(\Omega)} \| f^2 \|_{L^{\tilde q} (\Omega)}
= C \| g \|_{L^2(\Omega)}^{\frac{1}{q}}
\| f \|_{L^{\frac{4q}{2q-1}}(\Omega)}^2 ,
\end{eqnarray*}
where we have used
$\| \psi_{1z}\|_{L^{\tilde p}(\Omega)} \leq \| \psi_{1zz}\|_{L^{\tilde p}(\Omega)}$
which follows from the Poincare inequality along the $z$-direction,
and the weighted Calderon-Zygmund estimate
$\| \psi_{1zz} \|_{L^{\tilde p} (\Omega)} \leq C \| \omega_1 \|_{L^{\tilde p} (\Omega)}$,
which we have proved in Lemma 2. In order to apply Lemma 2, we need to
show that $\psi_1 \in L^{\tilde p}(\Omega)$. To see this, we use the
fact that
$V^r=-\frac{\partial}{\partial z} \psi^{\theta} \in L^m$ for any
$m\geq 2$ and $\psi^{\theta}$ is odd. Thus we can apply the Poincare
inequality along $z$-direction to show that $\psi^{\theta} \in L^m(\Omega)$,
which implies that $\psi_1 = \frac{1}{r} \psi^{\theta} \in L^s(\Omega)$ for
any $s>1$ in the region $r>1$. In the region where $r$ is bounded and
$0<z<2$, the interior estimate of the general elliptic theory implies that
$\psi_1 \in L^{\tilde p}(\Omega)$ whenever $\omega_1 \in L^{\tilde p}(\Omega)$.

Using the Sobolev interpolation inequality,
\begin{equation}
\|f \|_{L^{\frac{4q}{2q-1}}} \leq C \|f\|_{L^2}^\alpha
\| \nabla f \|_{L^2}^{1-\alpha} \;,
\label{Sob}
\end{equation}
with $\alpha = 1 - \frac{3}{4q}$ and $q \geq 1$, we obtain
\begin{equation}
\frac{1}{2p} \frac{d}{dt} \int f^2 r dr dz \leq C_p \|g\|_{L^2}^{\frac{1}{q}}
\|f\|_{L^2}^{2-\frac{3}{2q}}
\| \nabla f \|_{L^2}^{\frac{3}{2q}} \; - \;
\frac{(2p-1)}{p^2} \| \nabla f \|_{L^2}^2.
\label{eqn-2}
\end{equation}

By the assumption, the initial conditions for $\psi_1$,
$u_1$, and $\omega_1$ are odd functions in $z$. This
oddness property of the solution is preserved dynamically. Thus we have
$f(r,0,t) = g(r,0,t) = 0$. By using the Poincare inequality along the
$z$-direction, we obtain
\begin{equation}
\|f\|_{L^2(\Omega)} \leq  \| f_z \|_{L^2(\Omega)}, \quad
\|g\|_{L^2(\Omega)} \leq  \| g_z \|_{L^2(\Omega)}.
\label{Poin-3d}
\end{equation}
It follows from (\ref{Poin-3d}) and (\ref{eqn-2}) that
\begin{equation}
\frac{1}{2p} \frac{d}{dt} \int_{\Omega} f^2 r dr dz \leq C_p
\|g\|_{L^2(\Omega)}^{\frac{1}{q}}
\|\nabla f\|_{L^2(\Omega)}^2  \; - \;
\frac{(2p-1)}{p^2} \| \nabla f \|_{L^2(\Omega)}^2 \;.
\label{eqn-3}
\end{equation}

Next, we multiply $|\omega_1|^{2q-2}\omega_1$ to (\ref{w1-dyn}) and integrate over
$\Omega$. We have
\begin{eqnarray*}
\frac{1}{2q} \frac{d}{dt} \int_{\Omega} g^2 r dr dz &\leq &
\int_{\Omega} (u_1^2)_z |\omega_1|^{2q-2}\omega_1 r dr dz -
 \frac{(2q-1)}{q^2} \| \nabla g \|_{L^2(\Omega)}^2 \\
& \leq & \int_{\Omega} ( f^{\frac{2}{p}} )_z
|\omega_1|^{2q-2}\omega_1 r dr dz
- \frac{(2q-1)}{q^2} \| \nabla g \|_{L^2(\Omega)}^2 \\
&\leq  & (2-\frac{1}{q})
\int_{\Omega} f^{\frac{2}{p}} g^{1 - \frac{1}{q}} |g_z| r dr dz
-  \frac{(2q-1)}{q^2} \| \nabla g \|_{L^2(\Omega)}^2 \\
& \leq & (2-\frac{1}{q})
\| g_z \|_{L^2(\Omega)} \left ( \int_{\Omega} f^{\frac{4}{p}}
g^{2(1-\frac{1}{q})} r dr dz \right )^{1/2}
- \frac{(2q-1)}{q^2} \| \nabla g \|_{L^2(\Omega)}^2
 \;.
\end{eqnarray*}

Let $\tilde p = p/2$, $\tilde q = \frac{\tilde p}{\tilde p -1 } = \frac{p}{p-2}$.
Then we obtain by using the H\"older inequality that
\begin{eqnarray*}
\left ( \int_{\Omega} f^{\frac{4}{p}}
g^{2(1-\frac{1}{q})} r dr dz \right )^{1/2} & \leq &
\left ( \int_{\Omega} f^2 r dr dz \right )^{1/p}
\left ( \int_{\Omega} g^{2(1-\frac{1}{q})\frac{p}{p-2}} r dr dz \right )^{\frac{p-2}{2p}}\\
&\leq & \|f\|_{L^2(\Omega)}^{\frac{2}{p}} \|g \|_{L^2(\Omega)}^{1-\frac{2}{p}} \;,
\end{eqnarray*}
where we have used the relationship $p=2q$.
Therefore, we obtain
\begin{eqnarray}
\frac{1}{2q} \frac{d}{dt} \int_{\Omega} g^2 r dr dz &\leq &
(2-\frac{1}{q}) \| g_z\|_{L^2(\Omega)}
\| f \|_{L^2(\Omega)}^{\frac{2}{p}} \| g\|_{L^2(\Omega)}^{1-\frac{2}{p}}
- \frac{(2q-1)}{q^2} \| \nabla g \|_{L^2(\Omega)}^2 \nonumber \\
& \leq & C \| f \|_{L^2(\Omega)}^{\frac{4}{p}} \| g\|_{L^2(\Omega)}^{2(1-\frac{2}{p})}
-  \frac{(2q-1)}{2q^2} \| \nabla g \|_{L^2(\Omega)}^2 \nonumber \\
& \leq & C_q \frac{(2p-1)}{2p^2}\| f \|_{L^2(\Omega)}^2
+ \frac{(2q-1)}{4q^2} \| g \|_{L^2(\Omega)}^2 - \frac{(2q-1)}{2q^2}
\| \nabla g \|_{L^2(\Omega)}^2 \; \nonumber\\
& \leq & C_q \frac{(2p-1)}{2p^2} \| \nabla f \|_{L^2(\Omega)}^2
- \frac{(2q-1)}{4q^2} \| \nabla g \|_{L^2(\Omega)}^2 \; ,
\label{eqn-4}
\end{eqnarray}
where we have used the Poincare inequality (\ref{Poin-3d}) in the last inequality.

Multiplying (\ref{eqn-3}) by $C_q$ and adding
the resulting equation to (\ref{eqn-4}), we get
\begin{eqnarray}
&&\frac{d}{dt} \left (
C_q \int_{\Omega} f^2 r dr dz
+\frac{1}{2q} \int_{\Omega} g^2 r dr dz \right )\nonumber  \\
&& \leq
C_q \left (C_p \| g \|_{L^2(\Omega)}^{\frac{1}{q}}
\| \nabla f \|_{L^2(\Omega)}^2 -
\frac{(2p-1)}{p^2}\| \nabla f \|_{L^2(\Omega)}^2 \right )
 \nonumber\\
& & + C_q \frac{(2p-1)}{2p^2} \| \nabla f \|_{L^2(\Omega)}^2
- \frac{(2q-1)}{4q^2} \| \nabla g \|_{L^2(\Omega)}^2 \nonumber  \\
&& \leq C_q \left (C_p \| g \|_{L^2(\Omega)}^{\frac{1}{q}}
 - \frac{(2p-1)}{2p^2}  \right ) \| \nabla f \|_{L^2(\Omega)}^2 \;.
\label{eqn-5}
\end{eqnarray}
Using (\ref{UW-init1}), we have
\begin{eqnarray}
u_1(r,z,0) = \epsilon^{\delta} U_1 (\epsilon r,
z), \quad
\omega_1(r,z,0) = \epsilon^{\delta} W_1 (\epsilon r, z),
\end{eqnarray}
where $U_1(\tilde r ,\tilde z)$ and $W_1(\tilde r ,\tilde z)$ have
bounded $L^{2p}$ and $L^{2q}$ norms independent of $\epsilon$. We obtain
\begin{eqnarray}
\|g_0 \|_{L^2}^{\frac{1}{q}} & = & \epsilon^{\delta- \frac{1}{q}}
\| W_1 \|_{L^{2q}(\Omega)},\label{g0-bound}\\
\|f_0 \|_{L^2}^{\frac{1}{p}} & = & \epsilon^{\delta - \frac{1}{p}}
\| U_1 \|_{L^{2p}( \Omega)} \label{f0-bound}.
\end{eqnarray}

We would like to choose $\epsilon$ small enough to ensure that the right hand
side of (\ref{eqn-5}) is negative. To this end, we require that
\begin{equation}
C_p \|g\|_{L^2(\Omega)}^{\frac{1}{q}}(t) \leq
 \frac{(2p-1)}{2p^2}.
\label{cond1}
\end{equation}
If (\ref{cond1}) holds, we would get
\begin{eqnarray}
\frac{d}{dt} \left ( \frac{C_q}{2p}
\int_{\Omega} f^2 r dr dz +\frac{1}{2q} \int_{\Omega} g^2 r dr dz \right ) \leq 0 .
\end{eqnarray}
This implies that
\begin{eqnarray}
\frac{C_q}{2p} \int_{\Omega} f^2 r dr dz +
\frac{1}{2q} \int_{\Omega} g^2 r dr dz \leq
\frac{C_q}{2p} \int_{\Omega}
f_0^2 r dr dz +\frac{1}{2q} \int g_0^2 r dr dz .
\label{global-bound}
\end{eqnarray}
This gives a global bound on $\|u_1\|_{L^{2p}(\Omega)}$ and
$\|\omega_1\|_{L^{2q}(\Omega)}$
in terms of their initial data. In particular, we have by using
(\ref{g0-bound})-(\ref{f0-bound}) that
\begin{eqnarray}
\| g \|_{L^2(\Omega)}^{\frac{1}{q}} (t)
&\leq & \| g_0 \|_{L^2(\Omega)}^{\frac{1}{q}} +
(C_q/2)^{\frac{1}{2q}} \| f_0 \|_{L^2(\Omega)}^{\frac{2}{p}}
\nonumber \\ \nonumber \\
& \leq &
\epsilon^{\delta- \frac{1}{q}}
\| W_1 \|_{L^{2q}(\Omega)}
+ (C_q/2)^{\frac{1}{2q}} \epsilon^{2\delta - \frac{1}{q}}
\| U_1 \|_{L^{2p}( \Omega)}^2.
\label{bound-g}
\end{eqnarray}
Using (\ref{bound-g}), we conclude that the condition
(\ref{cond1}) is satisfied if we choose $\epsilon$ to satisfy
\begin{equation}
C_p \left (
\epsilon^{\delta- \frac{1}{q}}
\| W_1 \|_{L^{2q}(\Omega)} + (C_q/2)^{\frac{1}{2q}} \epsilon^{2\delta - \frac{1}{q}}
\| U_1 \|_{L^{2p}( \Omega)}^2 \right )
\leq  \frac{(2p-1)}{2p^2},
\label{cond2}
\end{equation}
which is true by taking $q > 1/\delta$ and $\epsilon$ small enough.
This condition would guarantee that the
global estimate (\ref{global-bound}) is valid for all time.
The global estimate (\ref{global-bound}) gives a global bound
on $\|u_1\|_{L^{2p}( \Omega)}$ and $\|\omega_1\|_{L^{2q}( \Omega)}$.
Using this global estimate, we can easily obtain a global bound on
$\|u^\theta\|_{L^{2p}}$ and $\|\omega^\theta\|_{L^{2q}}$ over a
bounded domain with $r\leq R$, which proves the global
regularity of the 3D Navier-Stokes equations (see Remark 3).
This completes the proof of Theorem 1.

\vspace{0.2in}
\noindent
{\bf Remark 4.}
The global regularity of the Navier-Stokes equations with initial
data (\ref{Data-3}) follows almost exactly as the proof of Theorem 1
except that the bounds for $\|g_0\|_{L^2}$ and $\|f_0\|_{L^2}$
are different. By using the same rescaling, we just need to consider the
rescaled initial data of the form
\begin{eqnarray}
u^\theta(r,z,0) = U_0 (\epsilon^\delta r,z), \quad
\omega^\theta(r,z,0) =  W_0 (\epsilon^\delta r, z).
\label{UW-init4}
\end{eqnarray}
Due to the difference in the scaling, the bounds for
$\|g_0\|_{L^{2q}}$ and $\|f_0\|_{L^{2q}}$ are now given by
\begin{eqnarray}
\|g_0 \|_{L^2}^{\frac{1}{q}} & = & \epsilon^{\delta (1-1/q)}
\| W_1 \|_{L^{2q}(\Omega)},\\
\|f_0 \|_{L^2}^{\frac{1}{p}} & = & \epsilon^{\delta (1-1/p)}
\| U_1 \|_{L^{2p}(\Omega)} .
\end{eqnarray}
Thus the condition (\ref{cond1}) is satisfied if we choose $q>1$ and $\epsilon$
small enough to satisfy
\begin{equation}
C_p \left (
\epsilon^{\delta (1- \frac{1}{q})}
\| W_1 \|_{L^{2q}(\Omega)} + (C_q/2)^{\frac{1}{p}} \epsilon^{\delta (2- \frac{1}{q})}
\| U_1 \|_{L^{2p}( \Omega)}^2 \right )
\leq  \frac{(2p-1)}{2p^2},
\label{cond3}
\end{equation}
Condition (\ref{cond3}) will guarantee that the
global estimate (\ref{global-bound}) is valid for all time.

\vspace{0.2in}
\noindent
{\bf Acknowledgments.} We would like to thank
Profs. Charles Fefferman and Fanghua Lin
for their interests in this work and for some stimulating
discussions. The research was in part supported by NSF
under the NSF FRG grant DMS-0353838 and ITR Grant ACI-0204932.

\bibliographystyle{amsplain}

\end{document}